\documentclass[12pt, a4paper, twoside]{article}

\usepackage[english]{babel}
\usepackage{amsmath}
\usepackage{latexsym}
\usepackage{amssymb,amsthm,amsfonts}
\usepackage[mathcal]{eucal}
\usepackage{amscd}
\usepackage[latin1]{inputenc}
\usepackage{color}
\def\cred{\color{red}}
\def\cblu{\color{blue}}

\let\cred\relax
\let\cblu\relax

\theoremstyle{plain}
\newtheorem{thm}{Theorem}[section]
\newtheorem{corol}[thm]{Corollary}
\newtheorem{lem}[thm]{Lemma}
\newtheorem{propos}[thm]{Proposition}
\theoremstyle{definition}

\theoremstyle{remark}
\newtheorem{remark}[thm]{Remark}

\newenvironment{theor}{\vspace{-6mm}\begin{itemize}\item[]\begin{thm}}{\end{thm}\end{itemize}}
\newenvironment{lemma}{\vspace{-6mm}\begin{itemize}\item[]\begin{lem}}{\end{lem}\end{itemize}}

\newenvironment{prop}{\vspace{-6mm}\begin{itemize}\item[]\begin{propos}}{\end{propos}\end{itemize}}

\newcommand{\R}{\mathbb{R}}        % per i reali
\newcommand{\vett}[2]{#1(0,T; #2)} % per gli spazi funzionali a valori vettoriali

\newcommand{\der}[2]{\frac{\displaystyle \partial #1}{\displaystyle \partial #2}} % per le derivate, in notazione leibnitziana
\newcommand{\abs}[1]{\left| #1 \right|}                                           % per il valore assoluto
\newcommand{\norm}[1]{\left\| #1 \right\|}                                        % per la norma
\newcommand{\scal}[2]{\left( #1 , #2 \right)}                                   % per il prodotto scalare
\newcommand{\dual}[2]{\left\langle #1 ,  #2 \right\rangle}                       % per la dualit?

\newcommand{\nh}[1]{\norm{#1}^2_H}                                                % per il quadrato della norma H
\newcommand{\Pb}{\left(\textbf{P}_{\beta}\right)}
\newcommand{\Pab}{\left(\textbf{P}_{\alpha, \hskip.5pt \beta}\right)}

\newcommand{\ya}{y_\alpha}
\newcommand{\ua}{u_\alpha}
\newcommand{\xa}{\xi_\alpha}
\newcommand{\yh}{\widehat{y}_\alpha}
\newcommand{\uh}{\widehat{u}_\alpha}
\newcommand{\xh}{\widehat{\xi}_\alpha}
\newcommand{\wza}{w_{0, \hskip.5pt \alpha}}
\newcommand{\uza}{u_{0, \hskip.5pt \alpha}}
\newcommand{\vza}{v_{0, \hskip.5pt\alpha}}
\newcommand{\wzh}{\widehat{w}_{0, \hskip.5pt\alpha}}
\newcommand{\uzh}{\widehat{u}_{0, \hskip.5pt\alpha}}
\newcommand{\vzh}{\widehat{v}_{0, \hskip.5pt\alpha}}
\newcommand{\fa}{f_\alpha}
\newcommand{\fh}{\widehat{f}_{\alpha}}

\renewcommand\div{\mathop{\rm div}\nolimits}
\numberwithin{equation}{section}

\begin{document}
\renewcommand{\baselinestretch}{1.05}

\pagestyle{myheadings}
\newcommand\testopari{\sc Giacomo Canevari and Pierluigi Colli}
\newcommand\testodispari{\sc Convergence properties for a phase field system}
\markboth{\testodispari}{\testopari}
\renewcommand{\thefootnote}{\fnsymbol{footnote}}

\thispagestyle{empty}
%\title{}

\begin{center}
{{\bf\huge Convergence properties\\[1.5mm]
for a generalization of the Caginalp\\[3mm]
phase field system\footnote{{\bf Acknowledgment.}\quad\rm
The financial support of the MIUR-PRIN Grant 2008ZKHAHN 
\emph{``Phase transitions, hysteresis 
and multiscaling''} and of the IMATI of CNR in
Pavia is gratefully acknowledged.}}}

\vspace{6mm}
{\large\sc Giacomo Canevari and Pierluigi Colli}

\vspace{3mm}
{\sl Dipartimento di Matematica  ``F. Casorati'', Universit\`a di Pavia\\%[1mm]
Via Ferrata, 1, 27100 Pavia,  Italy\\%[1mm]
E-mail: {\tt giacomo.canevari@gmail.com \ \  pierluigi.colli@unipv.it}}

\end{center}

\vskip6mm
\begin{abstract}

%\noindent {\bf Abstract.} 
We are concerned with a phase field system consisting of   
two partial differential equations in terms of the variables 
\emph{thermal displacement}, that is basically the time integration 
of temperature, and \emph{phase parameter}.  The system is a generalization of the
well-known Caginalp model for phase transitions, when including a diffusive term for 
the thermal displacement in the 
balance equation and when dealing with an arbitrary maximal monotone graph, 
along with a smooth anti-monotone function, in the phase equation. 
A Cauchy-Neumann problem has been studied for such a system in \cite{cc1}, by proving well-posedness and regularity results, as well as convergence of the problem 
as the coefficient of the diffusive term for the thermal displacement tends to zero.
The aim of this contribution is rather to investigate the asymptotic behaviour of the 
problem as the coefficient in front of the Laplacian of the temperature goes to $0$: 
this analysis is motivated by the types III and II cases in the thermomechanical theory 
of Green and Naghdi.  
Under minimal assumptions on the data of the problems, we show a convergence result. 
Then, with the help of uniform regularity estimates, we discuss the rate of convergence 
for the difference of the solutions in suitable norms. 
 
\end{abstract}

\section{Introduction} \label{intro} 
In this paper we consider the initial and boundary value problem 
% $\left(\textbf{P}_{\alpha, \hskip.5pt \beta}\right)$
\begin{equation}
w_{tt} - \alpha\Delta w_t - \beta\Delta w + u_t = f \quad \textrm{ in } \Omega
\times (0,T)
\label{1}
\end{equation}
\begin{equation}
u_t - \Delta u + \gamma (u)  + \overline g(u) \ni w_t \quad \textrm{ in } \Omega
\times (0,T)
\label{2}
\end{equation}
\begin{equation}
\partial_n w = \partial_n u = 0 \qquad \textrm{on } \Gamma\times (0, T)
\label{3}
\end{equation}
\begin{equation}
w(\hskip1pt\cdot \hskip1pt ,  0) = w_0\, , \quad w_t (\hskip1pt\cdot\hskip1pt,  0) = v_0 \, , \quad u(\hskip1pt\cdot\hskip1pt,  0) = u_0 \qquad \textrm{ in } \Omega
\label{4}
\end{equation}
where $\Omega \subset \R^3 $ is a bounded domain with smooth boundary $\Gamma$,
$T> 0$ represents some finite time, and $\partial_n $ denotes the outward 
normal derivative on $\Gamma$. Moreover, $ f$ is a given source term in equation~\eqref{1}, $\overline g : \R \to \R$ is a Lipschitz-continuous function, and  
$\gamma : \R \to 2^{\R} $ is a maximal monotone graph (cf., e.g.,  \cite{Brezis} or \cite{Barbu}): as $\gamma$ can be multivalued, in \eqref{2} there is the inclusion instead of the equality. Finally, $w_0,\,  v_0 , \, u_0$ stand for initial data.  

Problem~\eqref{1}--\eqref{4} 
% $\left(\textbf{P}_{\alpha, \hskip.5pt \beta}\right)$ 
has been studied in the paper \cite{cc1} 
when $\alpha$ and $\beta$ are two positive coefficients; also, still in \cite{cc1} 
the asymptotic behaviour of the problem as $\beta \searrow 0 $ has been investigated. 
The interest of this new contribution is to examine the asymptotic properties of %$\left(\textbf{P}_{\alpha, \hskip.5pt \beta}\right)$ 
\eqref{1}--\eqref{4} as $\alpha \searrow 0$, being $\beta > 0$ fixed once 
and for all. Indeed, both asymptotic investigations are interesting and deserve to be investigated, in such a way that 
the results proved in this paper give a completion to those of \cite{cc1}.

Let us point out that the coupled equations \eqref{1}--\eqref{2} form 
a system of phase field type. From the seminal work of Caginalp and coworkers (see, e.g., \cite{caginalp, cagnish}) it became clear that phase field systems are particularly suited to 
represent the dynamics of moving interfaces arising in thermally induced 
phase transitions. In our case, we consider the following expression for the 
total free energy
\begin{equation}
\Psi (\theta, u) = \int_\Omega \left( - \frac12 \theta^2 - \theta u + \phi (u) + G(u) + \frac12 |\nabla u |^2 \right)
\label{5}
\end{equation}
where the variables  $\theta$ and $u$ denote the (relative) temperature and the order parameter, respectively. We remark at once that our $w$ represents the thermal displacement variable, related to $\theta$~by 
\begin{equation}
w(\cdot, \, t) = w_0 + (1* \theta) (\cdot, \, t) = w_0 + \int_0^t\!\! \theta (\cdot, \, s)\, 
ds , \quad \ t\in [0,T].
\label{6}
\end{equation}
The integrand in \eqref{5} is characterized by the presence of different terms: the 
first one yields the concave purely caloric contribution, the second is the only term 
coupling the two variables, $\phi+ G$ is a function acting on $u$ and the last one is 
accounting for surface effects. The convex and lower semicontinuous function $ \phi : [0,+\infty] \to \R $ satisfies $\phi(0) = 0 = \min \phi$ and its subdifferential 
$\partial \phi$ coincides with $\gamma $, while $G$ stands for a smooth, 
in general concave, function such that $G' = \overline g$. A typical example 
for $\phi (u) + G(u)$  is offered by the double obstacle case, in which 
\begin{equation}
\phi(u) =  I_{[-1, +1]} (u) =  
\begin{cases} 0 & \text{if $|u|\leq 1$}
\\
+\infty &\text{if $|u| >1$}
\end{cases}
, \quad \ 
G(u) = 1- u^2,
\label{7}
\end{equation}
the two wells of $\phi + G $ being located in $ -1$ and $+1$: actually, 
one of the two is preferred as minimum of the potential in \eqref{5}, according to 
whether the temperature $\theta$ is negative or positive. The above example perfectly fits with the systematic view and designs of Michel Fr\'emond in non-smooth Thermomechanics~\cite{fremond}. In the case of \eqref{7}, the subdifferential
of the indicator function of the interval $[-1, +1]$ reads 
$$ \xi \in \partial I_{[-1, +1]} (u) \quad \hbox{ if and only if } \quad 
\xi \ \left\{
\begin{array}{ll}
\displaystyle
\leq \, 0 \   &\hbox{if } \ u=-1   
\\[0.1cm]
= \, 0 \   &\hbox{if } \ |u| < 1  
\\[0.1cm]
\geq \, 0 \  &\hbox{if } \  u = + 1  
\\[0.1cm]
\end{array}
\right. 
$$
so that here $\gamma = \partial I_{[-1, +1]}$ is really a graph with vertical lines. Let us also note that even a multi-well potential may be considered in \eqref{5}: indeed, 
it turns out that the function $ u \mapsto - \theta u + \phi (u) + G(u)$ can exhibit different wells according especially to the shape of $G$ in the region where 
$\phi \not= +\infty$. 

A different word used for the variable $w$ in \eqref{6} is ``freezing index'':
this terminology is motivated by earlier studies on the Stefan 
problem~\cite{duvaut, af} and preferred by some authors (cf.~\cite{magenes, fremond}).
We notice that a meaningful variable of the Stefan problem 
is the enthalpy $e$, which is defined as minus the 
variational derivative of $\Psi $ with respect to~$ \theta$:
\begin{equation}
\label{def-e}
e= - d_\theta \Psi  ,  \quad \hbox{ whence } \,  e = \theta + u = w_t + u . 
\end{equation}  
Then, the governing balance and phase equations are given~by 
\begin{equation}
e_{t} + \div {\bf q}  = f 
\label{1phys}
\end{equation}
\begin{equation}
u_t +  d_u \Psi =0
\label{2phys}
\end{equation}
where ${\bf q} $ denotes the thermal flux vector and $d_u \Psi$ stands for the 
variational derivative of $\Psi$ with  respect to $u$. Hence, \eqref{2phys} reduces 
exactly to \eqref{2} along with the Neumann homogeneous boundary condition for $u$. 
On the other hand, if we assume the classical Fourier law 
$$ {\bf q} = - \nabla \theta = - \nabla w_t $$
(for the moment let us take the heat conductivity coefficient just equal to 1), 
then  \eqref{1phys} is nothing but the usual energy balance equation as in the 
Caginalp model~\cite{caginalp}. One obtains the same equation as in the weak formulation  of the Stefan problem, in which however \eqref{2} is replaced by the mere pointwise inclusion $ u\in \left( \partial I_{[-1, +1]}\right)^{-1} ( \theta ) $, that equivalently reads $   \partial I_{[-1, +1]} ( u) \ni  w_t$.

An alternative approach to the Fourier diffusive law consists in adopting the so-called 
Cattaneo-Maxwell law (see, e.g., \cite{cgg1, colrec, MQ1} and references therein) in  which a time derivative of the flux appears:
\begin{equation}
{\bf q} + \varepsilon  {\bf q}_t = - \nabla \theta , \quad \hbox{ for } \, \varepsilon  
> 0 \, \hbox{ small}.
\label{qeps}
\end{equation} 
The Cattaneo-Maxwell law leads to the following equation 
\[
\varepsilon  \theta_{tt} + \theta_t - \Delta \theta +  u_t + \varepsilon  u_{tt}  = 
f + \varepsilon f_t \quad \textrm{ in } \Omega
\times (0,T),
\label{1mv}
\]
which has been investigated in  \cite{MQ1, MQ5}. On the other hand, solving \eqref{qeps} with respect to ${\bf q} $ yields
$$
{\bf q} = {\bf q_0} + k_\varepsilon* \nabla \theta , \ \hbox{ where } \, (k_\varepsilon * \nabla \theta) 
(x,t)  := \int_0^t  \!\! k_\varepsilon(t-s) \nabla \theta (x,s) ds ;
$$
${\bf q_0} (x,t)$ is known and can be incorporated in the source term, while 
$k_\varepsilon (t) $ is a given kernel depending on $\varepsilon $. Then, from \eqref{1phys} and a prescription like $ {\bf q} = {\bf q_0} + k * \nabla \theta$
for some fixed kernel $k$,  
we obtain the balance equation for the standard phase field model 
with memory which has a hyperbolic character and has been extensively studied in 
\cite{cgg1, cgg2}.

Green and Naghdi (see \cite{gn1,gn2,gn3,gn4})  presented an alternative 
treatment for a thermomechanical theory of deformable media. 
This theory takes advantage of an  
entropy balance rather than the usual entropy inequality. 
If we restrict our attention to the heat conduction, these authors
proposed three different hypotheses, labeled as type I, type II and type III, 
respectively. In particular, when type I is linearized, we recover the classical 
Fourier law (in terms of thermal displacement)
\begin{equation}
{\bf q} = - \alpha \nabla w_t  , \quad \alpha >0 .  
\label{typeI}
\end{equation} 
Type II is characterized by 
\begin{equation}
{\bf q} = - \beta \nabla w , \quad \beta >0 .
\label{typeII}
\end{equation} 
Moreover, a linearized version of type III reads
\begin{equation}
{\bf q} = - \alpha \nabla w_t - \beta \nabla w   . 
\label{typeIII}
\end{equation} 
We point out that \eqref{typeII}--\eqref{typeIII} laws
have been recently discussed, applied and compared by 
Miranville and Quintanilla in \cite{MQ2, MQ3, MQ4} (there the reader can find a rich 
list of references as well). In particular, our equation \eqref{1} reflects the use of \eqref{typeIII}, along with \eqref{1phys} and \eqref{def-e}, in deriving it; further, a no flux boundary condition for $ {\bf q}$ yields $ \partial_n w = 0  $ in~\eqref{3}.

Here we are: the system \eqref{1}--\eqref{4} comes as a consequence of  \eqref{1phys}--\eqref{2phys} when \eqref{5} and \eqref{typeIII} are postulated. Existence, continuous dependence, and  
regularity of the solution to the initial-boundary value 
problem~\eqref{1}--\eqref{4} have been studied in \cite{cc1}.  Observe that $\gamma$ is an arbitrary maximal monotone graph, therefore $\gamma$ may be multivalued, singular and with bounded domain. To this concern,
in addition to the example~\eqref{7} we also mention the singular 
case investigated in \cite{MQ4} for 
a convex potential $\phi$ with the same interval $[-1,+1]$ as domain, but with 
$$\phi (u)= \kappa_1 \big( (1+u) \ln (1+u) + (1-u) \ln (1-u) \big), $$ 
and for $\overline g(u) = G'(u) = - 2 \kappa_0 u $, where the meaningful constants $\kappa_0, \, \kappa_1$ have to satisfy $0<\kappa_1 <\kappa_0 $.

An interesting issue is of course the investigation of the asymptotic behaviour of the solutions to problem~\eqref{1}--\eqref{4} as one of the two parameters $\alpha$ 
or $\beta$ tends to $0$. In \cite{cc1} we have examined the limiting case $\beta 
\searrow 0$, keeping $\alpha>0 $ fixed 
and obtaining convergence of solutions to the analogous problem with 
$\beta=0$, which corresponds to \eqref{typeI}, i.e., the type~I case of Green and 
Naghdi. The paper \cite{cc1} also deals with
two error estimates of the difference of solutions in suitable norms, showing a  
rate of convergence that is linear with respect to $\beta$ in both estimates. 

When discussing the limit with $\alpha>0$  fixed, one can always exploit the properties of the parabolic term $- \alpha\Delta w_t$. This is no longer possible in the study of the limit as $\alpha \searrow 0$, which is the objective of the present paper. Indeed, if we let $\alpha \searrow 0$ and hold $\beta >0 $ fixed, 
we are led to the analog of problem~\eqref{1}--\eqref{4} with $\alpha =0$, in which 
equation \eqref{1} becomes hyperbolic. This limiting situation gives account of the type II hypothesis \eqref{typeII} of Green and Naghdi.

A change of variable turns out to be useful in obtaining estimates which does not depend on $\alpha$. Let $y$ denote the time-integrated enthalpy (cf.~\eqref{def-e}) 
$$y = w + 1*u.$$ 
Then, problem~\eqref{1}--\eqref{4} 
% $\left(\textbf{P}_{\alpha, \hskip.5pt \beta}\right)$ 
appears to be equivalent, in terms of $(y, \hskip.5pt u)$, to the problem
\begin{equation}\label{8}
y_{tt} -\alpha\Delta y_t - \beta\Delta y = -\alpha\Delta u - \beta\Delta(1*u) + f \quad \textrm{ in } \Omega\times (0,  T)
\end{equation}
\begin{equation}\label{9}
u_t - \Delta u + \gamma(u) + g(u) \ni y_t \quad \textrm{ in } \Omega\times (0,  T)
\end{equation}
\begin{equation}\label{10}
\partial_n y = \partial_n u = 0 \qquad \textrm{on } \Gamma\times (0,  T)
\end{equation}
\begin{equation}\label{11}
y(\hskip1pt\cdot\hskip1pt,  0) = w_0\, , \quad y_t (\hskip1pt\cdot\hskip1pt,  0) = v_0 + u_0 \, , \quad u(\hskip1pt\cdot\hskip1pt,  0) = u_0 \qquad \textrm{ in } \Omega \, ,
\end{equation}
where the function  $g$ is defined by $g(s) = \overline g(s) + s$ for all $s\in\R$ and results Lipschitz-continuous whenever $\overline g$ is. Let us notice that 
in the case $\alpha = 0$  problem~\eqref{8}--\eqref{11} reduces to a system for which well-posedness and regularity results have been proved in~\cite{cgg1, cgg2}.

By this change of variables, and dealing with some technicalities, we are able to provide a priori estimates on the solution of~\eqref{8}--\eqref{11}, independent of $\alpha$, which allow us to state a convergence result as $\alpha \searrow 0$ and prove it via compactness arguments. Then we investigate the bahaviour of the difference of the solutions to the problems with $\alpha >0 $ and $\alpha =0 $:
a first estimate on the convergence rate is obtained, although it is not linear but of order $1/2$ in $\alpha$. Stronger regularity estimates independent of $\alpha$ are then provided, and they enable us to produce also a linear estimate for the convergence rate, but only in the case when $\gamma$ is a smooth single-valued function defined on the whole real line. This assumption, although limitative, is not unrealistic, because it covers the physically relevant case of  the well-known `double-well' potential~\cite{caginalp} given by 
\begin{equation}
\phi (u) + G(u) = \frac\kappa 4 (u^2 -1)^2 \, , \quad \hbox{so that } \ 
\gamma (u) = \kappa \/u^3 \   \hbox{ and } \ g(u) = (1-\kappa)\/u \, ,
\label{double well}
\end{equation}
for some positive coefficient $\kappa$. 

The paper is organized as follows. In Section~\ref{no-ma} we present the main results related to the problem: the existence and uniqueness of a weak solution
for both cases  $\alpha >0 $ and $\alpha =0 $, the convergence theorem as $\alpha \searrow 0$, the regularity results yielding a strong solution and the first error estimate, the further regularity properties on 
the phase variable $u$ and the enthalpy $y_t$, and finally the second error estimate. Next, in Section~\ref{conv} we recall a useful lemma for parabolic problems and then prove the convergence result. Strong solutions and the first estimate on the order of convergence are discussed in Section~\ref{sec: strong}. 
The concluding Section~\ref{sm-po} is devoted to the case of a smooth potential $\phi$ and it contains the proofs of the additional regularity estimates and of the linear rate of convergence for the difference of solutions in suitable norms.

\section{Notation and main results}
\label{no-ma}

Before stating clearly the formulation of the problem and the main results we achieve, we recall some notation. Let $\Omega\subseteq\R^3$ be a bounded smooth domain, with boundary $\Gamma = \partial\Omega$, and let $T > 0$ be some final time. {\cred Set}
\[
Q_t = \Omega\times (0,  t) \, , \qquad Q = Q_T \, ,
\]
\[
H = L^2(\Omega) \, , \qquad V = H^1(\Omega) \, , \qquad W = \left\{v\in H^2(\Omega) : \: \partial_n v = 0 \textrm{ a.e. on } \Gamma \right\} \, .
\]
We embed $H$ in $V'$, by means of the $L^2(\Omega)$ inner product:
\[
\dual{h}{v} = \scal{h}{v}_H \qquad \textrm{for all } h\in H \, , \: v\in V \, ;
\]
we warn that $\norm{\cdot}_H$ means the norm of both $H$ and $H^3$. If $a$ is a function of space and time variables, we define the function $1*a$ by the formula
\[
(1 * a)(t) = \int_0^t a(s)\, ds \, , \qquad 0\leq t \leq T \, . 
\]
The symbol $c$ stands for different positive constants, depending on $\Omega$, $T$ and the data, which stay bounded whenever $\alpha, \, \beta$ vary in a bounded subset of $(0,  +\infty)$. {\cred In particular, for the sake of simplicity let us assume that
\begin{equation}
 0 < \alpha \leq 1 .\label{pc2}
\end{equation}
}% fine \cred
Any positive constant depending on the $\Omega$, $T$ and the data, which may tend to $+\infty$ when $\beta\searrow 0$ but stays bounded as $\alpha\searrow 0$, will be denoted by $c_\beta$.

For the sake of clarity, we call $\mathcal D_\alpha = \left\{\wza, \hskip.5pt \vza, \hskip.5pt \uza, \hskip.5pt \fa \right\}$, $(\ya, \hskip.5pt \ua, \hskip.5pt \xa)$ {\cred the family of data and the solution of the problem~\eqref{8}--\eqref{11} (cf. also \eqref{1}--\eqref{4}), 
while $\mathcal D = \left\{w_0, \hskip.5pt u_0, \hskip.5pt v_0, \hskip.5pt f\right\}$ and $(y, \hskip.5pt u, \hskip.5pt \xi)$ will denote the variables of related problem with} $\alpha = 0$.

Now we deal with the assumptions on data. We require that
\begin{equation}
0 < \alpha \leq 1\, , {\cblu \quad \beta > 0 }
\label{alpha, beta}
\end{equation}
\begin{equation}
\gamma\subseteq\R\times\R \: \textrm{ maximal monotone graph, with } \gamma(0) \ni 0
\label{gamma}
\end{equation}
\begin{equation}
\phi: \R\longrightarrow [0,  +\infty] \: \textrm{convex and lower semicontinuous, with } \phi(0) = 0 \textrm{ and } \partial\phi = \gamma
\label{phi}
\end{equation}
\begin{equation}
g: \R\longrightarrow\R \:\textrm{ Lipschitz-continuous}
\label{g}
\end{equation}
\begin{equation}
\fa \in \vett{W^{1,  1}}{V'} + \vett{L^1}{H}
\label{f}
\end{equation}
\begin{equation}
\wza \in V \, , \qquad \vza \in H \, , \qquad \uza \in H \, , \qquad \phi(\uza) \in L^1(\Omega) \, ;
\label{initial data}
\end{equation}
we also assume, for the norms of the data, a bound independent of $\alpha$:
\begin{equation}
\norm{\fa}_{\vett{W^{1,  1}}{V'} + \vett{L^1}{H}} + \norm{\wza}_V + \norm{\vza}_H + \norm{\uza}_H + \norm{\phi(\uza)}_{L^1(\Omega)} \leq c \, .
\label{d norm}
\end{equation}
We denote the effective domain of $\gamma$ with $D(\gamma)$. {\cred Now, we are ready to formulate precisely the problem as subsequent to the change of variable.}

\bigskip 
\noindent
\textbf{Problem $\Pab$.} Find $(\ya, \hskip.5pt \ua, \hskip.5pt \xa)$ satisfying
\begin{equation}
\ya \in \vett{W^{2,  1}}{V'} \cap \vett{W^{1,  \infty}}{H} \cap \vett{L^\infty}{V}
\label{y_alpha}
\end{equation}
\begin{equation}
\ua \in \vett{H^1}{V'} \cap C^0\left([0, T]; \, H \right) \cap \vett{L^2}{V}
\label{u_alpha}
\end{equation}
\begin{equation}
\xa\in L^2(Q) \, , \qquad \ua\in D(\gamma) \ \textrm{ and } \ \xa\in\gamma(\ua)\  \textrm{ a.e. in } \, Q
\label{xi_alpha}
\end{equation}
\begin{equation}
\begin{split}
\dual{\partial_t^2\ya(t)}{v} + \alpha\scal{\nabla \partial_t\ya(t)}{\nabla v}_H + \beta\scal{\nabla \ya(t)}{\nabla v}_H
= \alpha\scal{\nabla\ua (t)}{\nabla v}_H \\
+ \beta\scal{\nabla(1*\ua)(t)}{\nabla v}_H + \dual{\fa(t)}{v} \qquad\textrm{for all } v\in V \textrm{ and a.a. } t\in (0,  T)
\end{split}
\label{eq-A_alpha}
\end{equation}
\begin{equation}
\begin{split}
\dual{\partial_t\ua(t)}{v} + \scal{\nabla \ua(t)}{\nabla v}_H + \scal{\xa(t)}{v}_H + \scal{g(\ua)(t)}{v}_H
= \scal{\partial_t \ya(t)}{v}_H \\ \qquad \textrm{for all } v\in V \textrm{ and a.a. } t\in(0,  T)
\end{split}
\label{eq-B_alpha}
\end{equation}
\begin{equation}
\ya(0) = \wza\ \textrm{ in } \, V \, , \quad \partial_t\ya(0) = \vza + \uza \ \textrm{ in } \, V'\, , \quad \ua(0) = \uza \ \textrm{ in } \, H.
\label{i-c-alpha}
\end{equation}

We obtain the formulation of problem $\Pb$ by setting $\alpha = 0$ in problem $\Pab$.

\bigskip 
\noindent
\textbf{Problem $\Pb$.} Find $(y, \hskip.5pt u, \hskip.5pt \xi)$ satisfying
\begin{equation}
y \in \vett{W^{2,  1}}{V'} \cap \vett{W^{1,  \infty}}{H} \cap \vett{L^\infty}{V}
\label{y}
\end{equation}
\begin{equation}
u \in \vett{H^1}{V'} \cap C^0\left([0, T]; \, H \right) \cap \vett{L^2}{V}
\label{u}
\end{equation}
\begin{equation}
\xi\in L^2(Q) \, , \qquad u\in D(\gamma) \ \textrm{ and } \ \xi\in\gamma(u)\  \textrm{ a.e. in } \, Q
\label{xi}
\end{equation}
\begin{equation}
\begin{split}
\dual{\partial_t^2 y(t)}{v} + \beta\scal{\nabla y(t)}{\nabla v}_H = \beta\scal{\nabla(1*u)(t)}{\nabla v}_H + \dual{f(t)}{v} \\
\qquad\textrm{for all } v\in V \textrm{ and a.a. } t\in (0,  T)
\end{split}
\label{eq-A}
\end{equation}
\begin{equation}
\begin{split}
\dual{\partial_t u(t)}{v} + \scal{\nabla u(t)}{\nabla v}_H + \scal{\xi(t)}{v}_H + \scal{g(u)(t)}{v}_H = \scal{\partial_t y(t)}{v}_H \\
\qquad \textrm{for all } v\in V \textrm{ and a.a. } t\in(0,  T)
\end{split}
\label{eq-B}
\end{equation}
\begin{equation}
y(0) = w_0\ \textrm{ in } \, V \, , \quad \partial_t y(0) = v_0 + u_0 \ \textrm{ in } \, V'\, , \quad u(0) = u_0 \ \textrm{ in } \, H.
\label{i-c-}
\end{equation}

{\cred \begin{prop}[Well-posedness] \label{th: w-p}
Under the assumptions \eqref{alpha, beta}--\eqref{initial data} for both the families of data $\mathcal D_\alpha$ and $\mathcal D$, the two problems $\Pab$ and $\Pb$ admit unique solutions $(\ya, \hskip.5pt \ua, \hskip.5pt \xa)$ and $(y, \hskip.5pt u, \hskip.5pt \xi)$, respectively.
\end{prop}

For the well-posedness of problem $\Pab$ we refer to \cite[Theorem 2.1]{cc1}, while the proof of the corresponding result for {\cblu $\Pb$} can be found in \cite{cgg2},
with the warning that the reader should argue with $y - w_0$ (in place of $y$)
since in the framework of \cite{cgg2} the initial value for the related variable is $0$.

In this work, we are interested in the limit as $\alpha\searrow 0$ for problem $\Pab$; so, what really affect us are the estimates for $(\ya, \hskip.5pt \ua, \hskip.5pt \xa)$ independent of $\alpha$, rather then mere regularity results for problem $\Pb$. For this reason, all the results are stated in terms of the two families of data $\mathcal D_\alpha$ and $\mathcal D$.
} % \end of \cred

\begin{theor}[Convergence as $\alpha \searrow 0$] \label{th: convergence}
{\cred Assume \eqref{alpha, beta}--\eqref{d norm} and let the families of data $\mathcal D_\alpha$ and $\mathcal D$ satisfy}
\begin{equation}
\fa \rightharpoonup f \quad \textrm{in } \vett{W^{1,  1}}{V'} + \vett{L^1}{H}
\label{f convergence}
\end{equation}
\begin{equation}
\wza \rightharpoonup w_0 \quad \textrm{in } V \, , \qquad \vza \rightharpoonup v_0 \, , \: \uza \rightharpoonup u_0 \quad \textrm{in } H \, .
\label{d convergence}
\end{equation}
Then, the convergences
\[
\ya \rightharpoonup^* y \quad \textrm{in } \vett{W^{1,  \infty}}{H} \cap \vett{L^\infty}{V}
\]
\[
\ua \rightharpoonup u \quad \textrm{in } \vett{H^1}{V'} \cap 
\vett{L^2}{{\cred V}} \, , \qquad 
\xa \rightharpoonup \xi \quad \textrm{in } L^2(Q)
\]
hold {\cred as $\alpha\searrow 0$.}
\end{theor}

Now, we strenghten the hypotheses on the initial data and on $f$, in order to obtain a strong solution to problem $\Pb$. In addition, we are able to provide an estimate on the convergence error, when $\alpha\searrow 0$.
{\cred Let us suppose that}
\begin{equation}
\fa \in \vett{W^{1,  1}}{H} + \vett{L^1}{V}
\label{f reg1}
\end{equation}
\begin{equation}
\wza \in W \, , \qquad \vza \in V \, , \qquad \uza \in V 
\label{d reg1}
\end{equation}
\begin{equation}
\norm{\fa}_{\vett{W^{1,  1}}{H}+\vett{L^1}{V}} + \norm{\wza}_W + \norm{\vza}_V + \norm{\uza}_V \leq c \, .
\label{dn reg1}
\end{equation}

\begin{theor}[Regularity and strong solution] \label{th: reg1}
{\cred Assume \eqref{alpha, beta}--\eqref{d norm}.  In addition, let 
\eqref{f convergence}--\eqref{dn reg1} hold for $\mathcal D_\alpha$ and $\mathcal D$. 
Then} the solution 
$(y, \hskip.5pt u, \hskip.5pt \xi)$ of problem $\Pb$ fulfills
\begin{equation} \label{pc11}
y \in \vett{W^{2,  1}}{H} \cap \vett{W^{1,  \infty}}{V} \cap \vett{L^\infty}{W}
\end{equation}
\begin{equation}\label{pc12}
u \in \vett{H^1}{H} \cap C^0\left([0,  T]; \, V\right) \cap \vett{L^2}{W} \, .
\end{equation}
In particular, $(y, \hskip.5pt u, \hskip.5pt \xi)$ solves a strong reformulation of problem $\Pb$, namely,
\begin{equation}\label{pc13}
\partial_t^2 y - \beta \Delta y = -\beta \Delta(1 * u) + f \qquad \textrm{a.e. in } Q
\end{equation}
\begin{equation}\label{pc14}
\partial_t u - \Delta u + \xi + g(u) = \partial_t y \, , \quad \xi\in \gamma(u) \qquad \textrm{a.e. in } Q
\end{equation}
\begin{equation}\label{pc15}
\partial_n y = \partial_n u = 0 \qquad \textrm{a.e. on } \Gamma\times [0,  T] \, .
\end{equation}
\end{theor}

\begin{theor}[First error estimate] \label{th: error est1}
Let the assumptions \eqref{alpha, beta}--\eqref{d norm} and \eqref{f reg1}--\eqref{dn reg1} {\cred hold; if the inequalities
\begin{equation}
\norm{\fa - f}_{\vett{W^{1,  1}}{V'} + \vett{L^1}{H}} \leq c\,\alpha^{1/2}
\label{f rate est1}
\end{equation}
\begin{equation}
\norm{\wza - w_0}_V + \norm{\vza - v_0}_H + \norm{\uza - u_0}_H \leq c\,\alpha^{1/2}
\label{d rate est1}
\end{equation}
are fulfilled, then the following estimate
\[
\norm{\ya - y}_{\vett{W^{1,  \infty}}{H} \cap \vett{L^\infty}{V}} + \norm{\ua - u}_{\vett{L^\infty}{H} \cap \vett{L^2}{V}} \leq c_\beta\, \alpha^{1/2} \,
\]
holds for the convergence error.}
\end{theor}

{\cblu \begin{remark} Thanks to the informations we get on $\ua - u$, and {\cred being} $y = w + 1*u$, the estimate for $\ya - y$ {\cred entails} an analogous one for $w_\alpha - w$, which was the original variable of the problem. The same will happen for all the estimates presented in the paper.
\end{remark}}

The subsequent result provides an $L^\infty$ estimate on $u$. For {\cred $s\in D(\gamma)$, let $\gamma^0(s)$ be the element of $\gamma(s)$ having minimum} modulus; we ask that
\begin{equation}
\uza \in W \, , \qquad \uza \in D(\gamma) \quad \textrm{a.e. in } \Omega \, , \qquad \gamma^0(\uza)\in H
\label{d reg2}
\end{equation}
\begin{equation}
\norm{\uza}_W + \norm{\gamma^0(\uza)}_H \leq c \, .
\label{dn reg2}
\end{equation}

\begin{theor}[Regularity for $u$] \label{th: reg2}
{\cred Assume \eqref{alpha, beta}--\eqref{d norm} and let
$\mathcal D_\alpha$, $\mathcal D$ satisfy the conditions  \eqref{f convergence}--\eqref{dn reg1}, \eqref{d reg2}--\eqref{dn reg2}. 
Then the solutions  $(\ya, \hskip.5pt \ua, \hskip.5pt \xa)$, $(y, \hskip.5pt u, \hskip.5pt \xi)$ of the respective problems $\Pab,$ $\Pb$ fulfill
\begin{equation}
\ua, \, u \in \vett{W^{1,  \infty}}{H} \cap \vett{H^1}{V} \cap \vett{L^\infty}{W} \, .
\label{reg 2}
\end{equation}
}% fine \cred
\end{theor}

All the results stated above hold for a domain $\Omega\subseteq\R^N$, where the dimension $N$ is arbitrary. {\cblu Note in particular that \eqref{reg 2} implies 
$\ua, \, u \in C^0([0,  T]; \, H^s(\Omega))$ for all $s < 2$
(cf., e.g., \cite[Sect.~8, Cor.~4]{Simon}).
On the other hand, il we let $N \leq 3$, the Sobolev embedding $H^s(\Omega)\hookrightarrow C^0(\overline\Omega)$, holding 
for $s > 3/2$, ensures that}
\[
u \in C^0(\overline Q) \, .
\]

The subsequent results deal with the particular case, which is still physically significative, in which $\gamma$ is a smooth function, defined on the whole {\cred of $\R$. We require that}
\begin{equation}
\gamma: \: \R\longrightarrow \R \quad \textrm{is single-valued and locally Lipschitz-continuous}
\label{gamma smooth}
\end{equation}
\begin{equation}
\fa \in L^2(Q) \, , \qquad \uza \in H^3(\Omega) \cap W \, , \qquad \norm{\fa}_{L^2(Q)} + \norm{\uza}_{H^3(\Omega)} \leq c \, .
\label{d reg2.5}
\end{equation}
{\cred Let us point out that the conditions \eqref{gamma smooth} and \eqref{d reg2.5}
together entail \eqref{d reg2}--\eqref{dn reg2}, because $ \gamma^0 (\uza) 
= \gamma (\uza)$ remains bounded in $L^\infty (\Omega )\ (\hookrightarrow H) $.}

\begin{theor}[Further regularity on $u$] \label{th: reg2.5}
{\cred Assume \eqref{alpha, beta}--\eqref{d norm} and let
$\mathcal D_\alpha$ and $\mathcal D$ fulfill  \eqref{f convergence}--\eqref{dn reg1} and \eqref{d reg2.5}; we also ask that $\gamma$ satisfy \eqref{gamma smooth}. Then the solutions  
$(\ya, \hskip.5pt \ua, \hskip.5pt \xa)$ and $(y, \hskip.5pt u, \hskip.5pt \xi)$ of the respective problems $\Pab$ and $\Pb$ satisfy
\begin{equation}\nonumber
\ua, \, u \in \vett{H^2}{H} \cap 
\vett{W^{1,  \infty}}{V} \cap \vett{H^1}{W} \, .
\end{equation}
}% fine \cred
\end{theor}

Finally, we provide regularity enough to obtain an $L^\infty$ {\cred bound} 
on $\partial_t y$ and a {\cred better} estimate on the convergence error. Let
\begin{equation}
\fa \in \vett{W^{2,  1}}{H} + \vett{W^{1,  1}}{V}
\label{f reg3}
\end{equation}
\begin{equation}
\vza \in W \, , \qquad \alpha\Delta\vza + \beta\Delta\wza + \fa(0) \in V 
\label{d reg3}
\end{equation}
\begin{equation}
\norm{\fa}_{\vett{W^{2,  1}}{H} + \vett{W^{1,  1}}{V}} + \norm{\vza}_W + \norm{\alpha\Delta\vza + \beta\Delta\wza + \fa(0)}_V \leq c \, .
\label{dn reg3}
\end{equation}
\begin{theor}[Further regularity on {\cred $y$}] \label{th: reg3}
{\cred Let \eqref{alpha, beta}--\eqref{d norm}, \eqref{f convergence}--\eqref{dn reg1} and \eqref{gamma smooth}--\eqref{dn reg3} hold. Then the solutions  $(\ya, \hskip.5pt \ua, \hskip.5pt \xa)$ and 
$(y, \hskip.5pt u, \hskip.5pt \xi)$ of the respective problems 
$\Pab$ and $\Pb$ satisfy
\begin{equation}\label{reg 3}
\ya, \,  y \in \vett{W^{2,  \infty}}{V} \cap \vett{W^{1,\infty}}{W} \, .
\end{equation}
}% fine \cred
\end{theor}
\begin{theor}[Second error estimate] \label{th: error est2}
{\cred Assume  \eqref{alpha, beta}--\eqref{d norm}, \eqref{f reg1}--\eqref{dn reg1}, \eqref{gamma smooth}--\eqref{dn reg3} and let the inequalities
\begin{equation}
\norm{\fa - f}_{\vett{W^{1,  1}}{V'} + \vett{L^1}{H}} \leq c\,\alpha
\label{f rate est1-a}
\end{equation}
\begin{equation}
\norm{\wza - w_0}_V + \norm{\vza - v_0}_H + \norm{\uza - u_0}_V \leq c\,\alpha
\label{d rate est1-a}
\end{equation}
be satisfied. Then, the following estimate 
\begin{equation} \label{pc31}
\begin{split}
&\norm{\ya - y}_{\vett{W^{1,  \infty}}{H} \cap \vett{L^\infty}{V}} \\
&\qquad {}+ \norm{\ua - u}_{\vett{H^1}{H} \cap C^0\left([0,  T]; \, V\right) \cap \vett{L^2}{W}} \leq c_\beta\, \alpha
\end{split}
\end{equation}
holds true for the convergence error. Moreover, if 
\begin{equation}
\norm{\fa - f}_{\vett{W^{1,  1}}{H} + \vett{L^1}{V}} + \norm{\wza - w_0}_W 
+ \norm{\vza - v_0}_V \leq c\,\alpha^{1/2} \, ,
\label{d rate est2}
\end{equation}
then we also have 
\begin{equation}\label{pc32}
\norm{\ya - y}_{\vett{W^{1,  \infty}}{V} \cap \vett{L^\infty}{W}} \leq c_\beta\, \alpha^{1/2} \, .
\end{equation}
}% fine \cred
\end{theor}

\begin{remark} Focusing on regularity of the domain $\Omega$, all the proofs contained in this paper work if standards results on Sobolev embedding and elliptic regularity apply. For instance, if $\Omega\subseteq\R^3$ is a convex polyhedron, then all the 
results of this paper remain true.
\end{remark}

\section{Preliminary results and convergence as $\alpha\searrow 0$}
\label{conv}

For the {\cred reader's} convenience, let us recall some results widely exploited in this paper. If $a, \, b\in L^2(Q)$ and $\sigma > 0$ is an arbitrary parameter, {\cred by} the H\"older and Young inequalities it is easy to infer {\cred that}
\begin{equation}
\int_{Q_t} ab \leq \frac{1}{2\sigma} \int_0^t \nh {a(s)} ds + \frac{\sigma}{2} \int_0^t \nh{b(s)} ds \, .
\label{pc0}
\end{equation}
Furthermore, if $a\in\vett{H^1}{H}$, by the {\cblu fundamental theorem of calculus} and the H\"older inequality we have
\begin{equation}
\nh{a(t)} = \nh{a(0) + \int_0^t \partial_t a(s) ds} \leq 2 \nh{a(0)} + 2T \int_0^t \nh{\partial_t a(s)} ds \, .
\label{ftc}
\end{equation}
{\cblu Finally, the following lemma states a well-known property of parabolic problems. {\cred For the sake of completeness, let us report a short proof.} 
\begin{lemma} \label{lemma: parabolic}
Let $h\in L^2(Q)$ and $z_0\in V$. Then, the problem
\[
\partial_t z - \Delta z = h \ \textrm{ a.e. in } Q
\]
\[
\partial_n z = 0 \ \textrm{ a.e. on } \Gamma\times [0,  T] \, , 
\qquad z(0) = z_0 \ \textrm{ a.e. in } \Omega
\]
admits exactly one solution, which fulfills
\begin{equation}
\norm{z}_{\vett{H^1}{H} \cap C^0\left([0, T]; \, V\right) \cap \vett{L^2}{W}} \leq c \left\{ \norm{z_0}_V + \norm{h}_{L^2(Q)} \right\} \, .
\label{parabolic}
\end{equation}
\end{lemma}
\proof 
The existence of a solution for this problem can be proved via an approximation-
passage to the limit scheme, such as, for instance, the Faedo-Galerkin method. The 
estimates needed to pass to the limit, by compactness arguments, are now formally 
provided. We test the equation with $\partial_t z$ and integrate over 
$Q_t${\cred , obtaining}:
\[
\begin{split}
\int_0^t\nh{\partial_t z(s)}ds + \frac{1}{2} \nh{\nabla z(t)} \leq \frac{1}{2}\nh{\nabla z_0} + \frac{1}{2} \int_0^t \nh{h(s)}ds + \frac{1}{2}\int_0^t\nh{\partial_t z(s)}ds \, .
\end{split}
\]
Such an inequality, combined with \eqref{ftc}, directly yields
\[
\norm{z}_{\vett{H^1}{H} \cap \vett{L^\infty}{V}} \leq c\left\{ \norm{z_0}_V + \norm{h}_{L^2(Q)} \right\} \, .
\]
{\cred Then, by comparison in the equation and on account of the homogeneous boundary conditions}, we obtain the estimate on $\norm{z}_{\vett{L^2}{W}}$ as well.
Setting $f\equiv 0$ and $z_0 = 0$, the previous inequality imply $z\equiv 0$, thus {\cred entailing} the uniqueness of solutions.
\endproof}

\newcommand{\ye}{y_{\alpha, \hskip.5pt \varepsilon}}

Now we are ready to face the proof of Theorem \ref{th: convergence}. Firstly, we establish estimates independent of $\alpha$ on the solution $(\ya, \hskip.5pt \ua, \hskip.5pt \xa)$ of problem $\Pab$; then, we pass to the limit for a suitable subsequence, with compactness arguments. As we already know {\cred (cf.~\cite{cgg2})
that Problem $\Pb$ has a \emph{unique} solution}, we will obtain the 
convergence of the whole family $(\ya, \hskip.5pt \ua, \hskip.5pt \xa)$.

\bigskip
\textbf{First a priori estimate} 
{\cred We aim to test  \eqref{eq-A_alpha}
by $\partial_t\ya $ and \eqref{eq-B_alpha} by $\ua $: as the former could be only 
formal due to regularity \eqref{y_alpha}, we adopt the technique of {\cblu singular 
perturbations} to regularize $\ya $}; all the results we need can be found, for 
example, in \cite[Prop.~6.1 and 6.2, p.~505, Prop.~6.3, p. 506]{cgg2}. {\cred 
Hence, for $\varepsilon > 0$ and a.a. $ t\in (0,T)$}, we define $\ye(t)\in V$ 
as the solution of the elliptic problem
\[
\ye(t) + \varepsilon^2 J\ye(t) = \ya(t) \, ,
\]
where $J$ is the Riesz isomorphism $V \longrightarrow V'$. Since
\[
\ya\in\vett{W^{2,  1}}{V'} \cap \vett{W^{1,  \infty}}{H} \cap \vett{H^1}{V} \, ,
\]
{\cred then it follows that $\ye\in\vett{W^{2,  1}}{V}$ and moreover, as  $\varepsilon\searrow 0$,}
\begin{equation}
\ye \rightharpoonup^* \ya \ \textrm{in } \vett{W^{1,  \infty}}{H}
\label{prop. A}
\end{equation}
\begin{equation}
\ye \rightharpoonup \ya \ \textrm{ in } \vett{H^1}{V} \, , \quad \hbox{{\cred whence} }\ \ye \longrightarrow \ya \ \textrm{ in } C^0\left([0,  T]; \, V\right) \, .
\label{prop. B}
\end{equation}
We begin {\cred with} testing equation \eqref{eq-B_alpha} with $\ua(t) \in V$ and integrating over $Q_t$. {\cred Recalling} that $\gamma$ is monotone and $g$ is Lipschitz-continuous, it is straightforward to get
\begin{equation}
\begin{split}
&\frac{1}{2}\nh{\ua(t)} + \int_0^t\nh{\nabla \ua(s)} ds \\
&\leq \frac{1}{2}\nh{\uza} + \frac{1}{2}\int_0^t\nh{\partial_t \ya(s)} ds
+ c\int_0^t \nh{\ua(s)} ds + c \, . 
\end{split}
\label{ts-1, 1}
\end{equation}
Now we test equation \eqref{eq-A_alpha} with $\partial_t\ye(t)$, integrate over $Q_t$, and pass to the limit as $\varepsilon{\cred \searrow}0$; it is convenient to treat each term separately. {\cred Thus,} by \cite[Prop. 6.2]{cgg2} and {\cblu \eqref{prop. B}}, we have
\[
\lim_{\varepsilon{\cred \searrow}0} \int_0^t\dual{\partial_t^2\ya(s)}{\partial_t\ye(s)} ds = \frac{1}{2}\nh{\partial_t \ya(t)} - \frac{1}{2}\nh{\vza + \uza} \, ,
\]
\[
\lim_{\varepsilon{\cred \searrow}0} \alpha\int_{Q_t} \nabla \partial_t \ya \cdot \nabla\partial_t\ye = \alpha \int_0^t \nh{\nabla \partial_t \ya(s)} ds \, ,
\]
and, {\cred since $(v, z) \mapsto (\nabla v,  \nabla z)_H$ gives a bilinear, symmetric, continuous and weakly coercive form in $V \times V$}, by \cite[Prop. 6.3]{cgg2} we find {\cred out that}
\[
\lim_{\varepsilon{\cred \searrow}0} \beta\int_{Q_t} \nabla \ya \cdot \nabla\partial_t \ye = \frac{\beta}{2}\nh{\nabla \ya(t)} - \frac{\beta}{2}\nh{\nabla \wza} \, .
\]
{\cred Before dealing with the terms on the right-hand side, let us 
sum up and see that
\begin{equation}
\begin{split}
&\frac{1}{2}\nh{\partial_t \ya(t)} + \alpha 
\int_0^t \nh{\nabla \partial_t \ya(s)} ds + \frac{\beta}{2} \norm{\ya(t)}^2_V \\
&\leq \frac{1}{2}\nh{\vza + \uza} + \frac{\beta}{2}\nh{\nabla \wza} + \frac{\beta}{2} \norm{\ya(t)}^2_H \\
&+ \lim_{\varepsilon{\cred \searrow}0} \left\{
\alpha\int_{Q_t} \nabla \ua \cdot \nabla\partial_t \ye
+ \beta\int_{Q_t} \nabla(1 * \ua) \cdot \nabla\partial_t \ye
+ \int_0^t \dual{\fa (s)}{\partial_t \ye(s)}ds \right\} \, ,  
\end{split}
\label{pc1}
\end{equation}
where we have added the same quantity (cf.~\eqref{ftc})
$$ 
\frac{\beta}{2} \norm{\ya(t)}^2_H  \ \left( \leq \beta \nh{\wza} + T\beta 
\int_0^t \nh{\partial_t \ya(s)} ds \right) 
$$
to both sides.} We now deal with the right-hand side of {\cred \eqref{pc1} and
note that}
\[
%\begin{split}
\lim_{\varepsilon{\cred \searrow}0} \alpha\int_{Q_t} \nabla \ua \cdot \nabla\partial_t \ye = \alpha\int_{Q_t}\nabla \ua \, \nabla\partial_t \ya %\\
\leq \alpha\int_0^t \nh{\nabla \ua(s)} ds + \frac{\alpha}{4}\int_0^t \nh{\nabla \partial_t \ya(s)} ds \, .
%\end{split}
\]
The other term involving $\ua$ is treated by time integration by parts{\cred : 
thanks to \eqref{ftc} and \eqref{pc0} we have that}
\[
\begin{split}
\lim_{\varepsilon{\cred \searrow}0} \beta\int_{Q_t} \nabla(1 * \ua) \cdot \nabla\partial_t \ye = \beta \lim_{\varepsilon{\cred \searrow}0} \left\{ \scal{\nabla(1 * \ua)(t)}{\nabla \ye(t)}_H - \int_{Q_t} \nabla \ua \cdot 
\nabla\ye \right\} \\
\leq  5\beta T\int_0^t \nh{\nabla \ua(s)} ds + \frac{\beta}{8}\nh{\nabla \ya(t)} + \frac{\beta}{4T} \int_0^t \nh{\nabla \ya(s)} ds 
\end{split}
\]
Finally, we fix a decomposition
\begin{equation}
\label{pc4}
\fa = \fa^{(1)} + \fa^{(2)} \, , \quad {\cred \textrm{with}  \quad
\norm{\fa^{(1)}}_{  \vett{W^{1,  1}}{V'} } + \norm{\fa^{(2)}}_{\vett{L^1}{H} } \leq c \, .}
\end{equation}
{\cred
We integrate by parts in time the term containing $f^{(1)}$ and exploit 
\eqref{i-c-alpha} and  \eqref{prop. A}--\eqref{prop. B}
to obtain
\begin{equation}\label{term-f1}
\begin{split}
&\lim_{\varepsilon{\cred \searrow}0} \int_0^t \dual{\fa^{(1)}(s)}{\partial_t \ye(s)}ds \\
&= \lim_{\varepsilon{\cred \searrow}0} \left\{ \dual{\fa^{(1)}(t)}{\ye(t)} - \dual{\fa^{(1)}(0)}{\ye(0)} - \int_0^t \dual{\partial_t \fa^{(1)}(s)}{\ye(s)} ds
\right\} \\
&\leq \frac{2}{\beta}\norm{\fa^{(1)}(t)}_{V'}^2 + \frac{\beta}{8}\norm{\ya(t)}_V^2 
+ \frac{1}{2} \norm{\fa^{(1)}(0)}_{V'}^2 \\ 
&\qquad+ \frac{1}{2} \norm{\wza}_V^2 + \int_0^t\norm{\partial_t \fa^{(1)}(s)}_{V'}\norm{\ya(s)}_V ds \\
&\leq c_\beta \norm{\fa^{(1)}}_{\vett{L^\infty}{V'}}^2 
+ \frac{\beta}{8}\norm{\ya(t)}_V^2 + \frac{1}{2} \norm{\wza}^2_V + \int_0^t \norm{\partial_t \fa^{(1)}(s)}_{V'}\norm{\ya(t)}_V ds  \, .
\end{split}
\end{equation}
Recalling \eqref{prop. A}, we easily get
\[
\lim_{\varepsilon{\cred \searrow}0}  \int_0^t \dual{\fa^{(2)}(s)}{\partial_t \ye(s)}ds = \int_{Q_t} \fa^{(2)} \, \partial_t \ya \leq \int_0^t\norm{\fa^{(2)}(s)}_H \norm{\partial_t \ya(s)}_H ds \, .
\]
Now, we collect all the inequalities referring to \eqref{pc1}. Owing to 
\eqref{pc2}, \eqref{d norm} and \eqref{pc4}, we infer that
\begin{equation}
\begin{split}
&\frac{1}{2}\nh{\partial_t \ya(t)} + \frac{3\alpha}{4} \int_0^t \nh{\nabla \partial_t \ya(s)} ds + \frac{\beta}{4} \norm{ \ya(t)}_V^2 \\
&\leq  c_\beta  + c \int_0^t\nh{\partial_t \ya(s)}ds
+ (1+ 5\beta T)  \int_0^t \nh{\nabla \ua(s)}ds 
 + c \int_0^t \nh{\nabla \ya(s)}ds \\
&\qquad 
+ \int_0^t \norm{\partial_t \fa^{(1)}(s)}_{V'}\norm{\ya(t)}_V ds
+ \int_0^t\norm{\fa^{(2)}(s)}_H\norm{\partial_t \ya(s)}_H ds  \, .
\end{split}
\label{ts-1, 2}
\end{equation}
Then, let us multiply \eqref{ts-1, 1} by $ (2 + 5\beta T) $ and add the resulting inequality to \eqref{ts-1, 2}, obtaining
\[
\begin{split}
&\nh{\ua(t)} + \int_0^t\nh{\nabla \ua(s)} ds + \frac{1}{2}\nh{\partial_t \ya(t)} + \frac{\alpha}{2} \int_0^t \nh{\nabla \partial_t \ya(s)} ds 
+ \frac{\beta}{4} \norm{ \ya(t)}_V^2 \\
& \leq c_\beta\left\{ 1+  \int_0^t \nh{\ua(s)}ds + \int_0^t\nh{\partial_t \ya(s)}ds + \int_0^t \nh{\nabla \ya(s)}ds \right\} \\
&\qquad{}+ \int_0^t \norm{\partial_t \fa^{(1)}(s)}_{V'}\norm{\ya(t)}_V ds
+ \int_0^t\norm{\fa^{(2)}(s)}_H\norm{\partial_t \ya(s)}_H ds \, ;
\end{split}
\]
at this point, a generalised version of the Gronwall lemma, which is stated in \cite[Teorema~2.1, p.~245]{baiocchi}, allows us to conclude that}
\begin{equation}
\norm{\ya}_{\vett{W^{1,  \infty}}{H} \cap \vett{L^\infty}{V}} + \sqrt{\alpha}\norm{\ya}_{\vett{H^1}{V}} + \norm{\ua}_{\vett{L^\infty}{H} \cap \vett{L^2}{V}} \leq c_\beta \, .
\label{est.1, 1}
\end{equation}

{\cred
\begin{remark}
Estimate \eqref{est.1, 1} can be easily modified in order to provide the continuous dependence on the data for problem $\Pab$. Indeed, taking the difference between the equations associated with two sets of data $\mathcal D_{\alpha, \hskip.5pt 1}$ and $\mathcal D_{\alpha, \hskip.5pt 2}$, the resulting system is formally analougous to $\Pab$ and can be treated similarly; the only difference is the term in $g(\ua)$, which is replaced by $g(u_{\alpha, \hskip.5pt 1}) - g(u_{\alpha, \hskip.5pt 2})$. We obtain an estimate where the constant in the right-hand side is independent of $\alpha$.
\end{remark}
}% fine \cred

\bigskip
\textbf{Second a priori estimate.} {\cred In order to prove the next estimate,} 
we need to approximate the graph $\gamma$ with its Yosida regularization; 
{\cred therefore}, for all $\varepsilon > 0$ we define
\[
\gamma_\varepsilon := \frac{1}{\varepsilon}\left\{I - \left(I + \varepsilon\gamma\right)^{-1}\right\} \quad \hbox{ and } \quad
\phi_\varepsilon(s) := \min_{\tau\in\R}\left\{\frac{1}{2\varepsilon}\abs{\tau - s}^2 + \phi(\tau)\right\} \quad  \hbox{ for }\,  s\in\R \, ,
\]
where $I$ denotes the identity on $\R$. We remind that $\phi_\varepsilon$ is a nonnegative {\cred and differentiable function, $\gamma_\varepsilon$ is increasing  and Lipschitz-continuous: moreover, we have that
\begin{equation}
\phi_\varepsilon ' (s) = \gamma_\varepsilon (s)\, , \quad 0 \leq \gamma_\varepsilon'(s) \leq \frac1\varepsilon  \, , \quad 0 \leq \phi_\varepsilon (s) \leq \phi (s) \quad \textrm{ for all
 }  s\in\R \, \hbox{ and } \, \varepsilon > 0
\label{properties}
\end{equation}
}% fine \cred
(see, e.g., \cite[Prop.~2.6, p.~28 and Prop.~2.11, p.~39]{Brezis} or \cite[pp.~57--58]{Barbu}). We establish an estimate for the solution $(\ya, \hskip.5pt \ua, \hskip.5pt \xa)$ to problem $\Pab$, in which we have replaced $\gamma$ by $\gamma_\varepsilon$; then, the same estimate will hold for the original problem, because of the passage to the limit {\cred as $\varepsilon \searrow 0$ already detailed} in \cite{cc1}.

We notice that $\gamma_\varepsilon(u) (t) \in V$ for a.a. $t$, due to the Lipschitz-continuity of $\gamma_\varepsilon$, and so it is an admissible test function for equation \eqref{eq-B_alpha}. The Lipschitz-continuity of $g$ and the formula 
\[
\gamma_\varepsilon(\ua)\,\partial_t\ua = \der{}{t}\phi_\varepsilon(\ua)
\]
allow us to write {\cred
\[
\begin{split}
&\norm{\phi_\varepsilon(\ua)(t)}_{L^1(\Omega)} + \int_0^t \gamma'_\varepsilon(\ua)\abs{\nabla \ua}^2 + \int_0^t\nh{\gamma_\varepsilon(\ua) (s)}ds \\
&\leq \norm{\phi_\varepsilon(\uza)}_{L^1(\Omega)} 
+ \int_0^t \nh{\partial_t\ya(s)} ds 
+ c\int_0^t\nh{\ua(s)} ds  +  c
+ \frac{1}{2}\int_0^t\nh{\gamma_\varepsilon(\ua) (s)}ds \, .
\end{split}
\]
}% fine \cred
We ignore the second term in the left-hand side, which is positive 
{\cred because of $\gamma'_\varepsilon \geq 0$, and use \eqref{properties},  \eqref{d norm} and \eqref{est.1, 1} to get}
\[
\norm{\phi_\varepsilon(\ua)(t)}_{L^1(\Omega)} + \frac{1}{2} \int_0^t\nh{\gamma_\varepsilon(\ua) (s)}ds \leq c_\beta \, .
\]
{\cred As} $\phi_\varepsilon\longrightarrow \phi$ pointwise (see \cite[Prop.~2.11, p.~39]{Brezis}), {\cred by applying the Fatou lemma to the first term  and passing to the limit as $\varepsilon \searrow 0$} we obtain
\begin{equation}
\norm{\phi(\ua)}_{\vett{L^\infty}{L^1(\Omega)}} + \norm{\xa}_{L^2(Q)} \leq c_\beta \, .
\label{est.2, 1}
\end{equation}
By comparison in the equation \eqref{eq-B_alpha}, we also {\cred deduce that}
\begin{equation}
\norm{\partial\ua}_{\vett{L^2}{V'}} \leq c_\beta \, .
\label{est.2, 2}
\end{equation}

\bigskip
\textbf{Passage to the limit as $\alpha\searrow 0$.} According to the estimates \eqref{est.1, 1}, \eqref{est.2, 1} and \eqref{est.2, 2}, it is fair to assume, {\cred up to subsequences},
\[
\ya \rightharpoonup^* y \quad \textrm{in } \vett{W^{1,  \infty}}{H} \cap \vett{L^\infty}{V}
\]
\[
\ua \rightharpoonup u \quad \textrm{in } \vett{H^1}{V'} \cap \vett{L^2}{V} \, , \qquad \xa \rightharpoonup \xi \quad \textrm{in } L^2(Q) \, .
\]
The generalised Ascoli theorem and the Aubin-Lions lemma {\cred (see, e.g., \cite[Sect.~8, Cor.~4]{Simon})} thus imply
\[
\ya \longrightarrow y \quad \textrm{strongly in } C^0\left([0,  T]; \, H\right)
\]
\[
\ua \longrightarrow u \quad \textrm{strongly in } C^0\left([0,  T]; \, V'\right) \ \textrm{and in } L^2(Q) \, ,
\]
and the Lipschitz-continuity of $g$ yields
\[
g(\ua) \longrightarrow g(u) \quad \textrm{strongly in } L^2(Q) 
\]
{\cred as $\alpha \searrow 0$. With all this information and in view of \eqref{f convergence}, \eqref{d convergence} it is not difficult to check that the limits $y, \, u, \, \xi$ 
form a triplet solving problem $\Pb$. The proof is analogous to the one developed in \cite[Section~5]{cc1}: in particular, let us point out that 
\[
\lim_{\alpha\searrow 0}\int_Q \xa \ua = \int_Q \xi u 
\]
which entails (see, e.g., \cite[Prop.~1.1, p.~42]{Barbu}) \eqref{xi}. Moreover, the regularity $\vett{W^{2, 1}}{V'}$ for $y$ can be recovered a posteriori by comparing terms of \eqref{eq-A}.}

\section{Strong solutions and first error estimate}
\label{sec: strong}

{\cred This section is devoted to the proofs of Theorems~\ref{th: reg1} 
and~\ref{th: error est1}. Regularity results for problem $\Pb$} are obtained by establishing estimates for $(\ya, \hskip.5pt \ua, \hskip.5pt \xa)$ independent of $\alpha$; then {\cred Theorem~\ref{th: convergence} and the weak or weak star compactness yield} the desired regularity for $(y, u, \xi)$. {\cred About the end 
of the section, we also discuss briefly how the thesis of Theorem~\ref{th: reg2}
follows from the furher requirements in \eqref{d reg2}--\eqref{dn reg2}.} 

We assume that $\mathcal D_\alpha$ and $\mathcal D$ satisfy hypotheses 
\eqref{f reg1}--\eqref{dn reg1}. Under these assumptions, {\cred from 
\cite[Theorem~2.2]{cc1} we know that
\begin{equation}
\ya \in \vett{W^{2,1}}{H} \cap \vett{W^{1, \infty}}{V} \cap \vett{H^1}{W}
%\nonumber
\label{pc5}
\end{equation}
\begin{equation}
\ua \in \vett{H^1}{H} \cap C^0\left([0, T]; \, V \right) \cap \vett{L^2}{W}
%\nonumber
\label{pc6}
\end{equation}
and $(\ya, \hskip.5pt \ua, \hskip.5pt \xa)$ solves Problem~$\Pab$ in a strong sense, that is, $\ya$ and $\ua$ satisfy
\begin{equation}\label{pc7}
\partial_t^2 \ya - \alpha\Delta \partial_t \ya - \beta\Delta \ya = -\alpha\Delta \ua -\beta\Delta(1 * \ua) + \fa  \quad \textrm{ a.e. in } Q
\end{equation}
\begin{equation}\label{pc8}
\partial_t \ua - \Delta \ua +  \xa + g(\ua ) = \partial_t \ya , \quad {\xa\in \gamma(\ua)}  \quad \textrm{ a.e. in } Q
\end{equation}
\begin{equation}\label{pc9}
\partial_n \ya = \partial_n \ua = 0 \quad \textrm{ a.e. on } \Gamma\times (0, T) 
\end{equation}
along with \eqref{i-c-alpha}, for every $\alpha \in (0,1)$.

\bigskip
\textbf{Third a priori estimate.} Thanks to \eqref{pc8}--\eqref{pc9}, we can apply 
Lemma~\ref{lemma: parabolic} with $z=\ua,$  $z_0 =\uza$, and $h=   \partial_t \ya
-  \xa - g(\ua ) $. Indeed, in view of the condition in {\cblu \eqref{dn reg1}} for the initial datum, and owing to \eqref{g} and to estimates \eqref{est.1, 1} and \eqref{est.2, 1}, from \eqref{parabolic} it follows that 
\begin{equation}
\norm{\ua}_{\vett{H^1}{H} \cap \vett{L^\infty}{V} \cap \vett{L^2}{W}} 
\leq c_\beta \, .
\label{est.3, 2}
\end{equation}
}% fine \cred

\bigskip
\textbf{Fourth a priori estimate.}
{\cred By virtue of \eqref{pc6}, we can plainly multiply equation \eqref{pc8} by $-\Delta\partial_t\ya$ and integrate over $Q_t$. Thus,} we obtain
\begin{equation}
\begin{split}
\frac{1}{2}\nh{\nabla\partial_t\ya(t)} + 
\alpha\int_0^t\nh{\Delta\partial_t\ya(s)}ds + \frac{\beta}{2}\nh{\Delta\ya(t)} \leq 
\frac{1}{2}\nh{\nabla (\vza + \uza) } \\
{}+ \frac{\beta}{2}\nh{\Delta\wza} + \alpha\int_{Q_t} \Delta\ua \, 
\Delta\partial_t\ya + \beta\int_{Q_t}\Delta(1*\ua) \, \Delta\partial_t\ya  - 
\int_{Q_t} \fa \, \Delta\partial_t\ya \, .
\end{split}
\label{ts-3, 1}
\end{equation}
A simple application of the H\"older inequality allow us to {\cred infer that}
\[
\begin{split}
\alpha\int_{Q_t} \Delta\ua \, \Delta\partial_t\ya \leq \frac{\alpha}{2} \int_0^t\nh{\Delta\ua(s)}ds + \frac{\alpha}{2} \int_0^t\nh{\Delta\partial_t\ya(s)}ds \, .
\end{split}
\]
We deal with the subsequent term integrating by parts in time:
\[
\begin{split}
\beta\int_{Q_t}\Delta(1*\ua) \, \Delta\partial_t\ya \leq 2\beta\nh{1*\Delta\ua(t)} + \frac{\beta}{8}\nh{\Delta\ya(t)} - \beta\int_{Q_t}\Delta\ua \, \Delta\ya \\
\leq {\cred \beta \left(4 T + 1\right) }\int_0^t \nh{\Delta\ua(s)} ds + \frac{\beta}{8}\nh{\Delta\ya(t)} + {\cred \frac{\beta}{4} } \int_0^t\nh{\Delta\ya(s)} ds \, .
\end{split}
\]
Now, we choose a {\cred split}
\begin{equation} \label{pc10}
\fa = \fa^{(1)} + \fa^{(2)} \, , \ \quad \textrm{with } \quad {\cred \norm{\fa^{(1)}}_{ \vett{W^{1,  1}}{H}} + \norm{ \fa^{(2)}}_{\vett{L^1}{V}} \leq c} \, ,  
\end{equation}
in order to estimate the term {\cred with} $\fa$. Concerning $\fa^{(1)}$, we integrate by parts in time and recall \eqref{dn reg1}: %\eqref{}
\[
\begin{split}
- \int_{Q_t} \fa^{(1)} \, \Delta\partial_t\ya = - \scal{\fa^{(1)}(t)}{\Delta\ya(t)}_H + \scal{\fa^{(1)}(0)}{\Delta\wza}_H + \int_{Q_t} \partial_t\fa^{(1)}\, \Delta \ya \\
\leq \frac{2}{\beta} \nh{\fa^{(1)}(t)} + \frac{\beta}{8} \nh{\Delta\ya(t)} +   \int_0^t \norm{\partial_t\fa^{(1)}(s)}_H \norm{\Delta\ya(s)}_H  ds + c \, .
\end{split}
\]
{\cred An integration by parts in space leads to}
\[
- \int_{Q_t} \fa^{(2)} \, \Delta\partial_t\ya \leq \int_0^t \norm{\nabla\fa^{(2)}(s)}_H \norm{\nabla\partial_t\ya(s)}_H ds \, .
\]
Collecting all these {\cred inequalities and taking advantage of \eqref{dn reg1}, \eqref{pc2}, \eqref{est.3, 2} and  \eqref{pc10}, 
from \eqref{ts-3, 1} we derive
\begin{equation}\label{pc23}
\begin{split}
\frac{1}{2}\nh{\nabla\partial_t\ya(t)} + \frac{\alpha}{2}\int_0^t\nh{\Delta\partial_t\ya(s)}ds + \frac{\beta}{4}\nh{\Delta\ya(t)} 
\leq  c_\beta + c\int_0^t\nh{\Delta\ya(s)} ds \\
%&\qquad 
+ \int_0^t \norm{\partial_t\fa^{(1)}(s)}_H \norm{\Delta\ya(s)}_H  ds 
+ \int_0^t \norm{\nabla\fa^{(2)}(s)}_H \norm{\nabla\partial_t\ya(s)}_H ds  \, .
\end{split}
\end{equation}
}% fine \cred
Hence, by means of the generalised Gronwall lemma (see{\cred , e.g.,} \cite{baiocchi}) {\cred and regularity properties for elliptic problems}, we obtain
\begin{equation}
\norm{\ya}_{\vett{W^{1,  \infty}}{V} \cap \vett{L^\infty}{W}} + \sqrt{\alpha}\norm{\ya}_{\vett{H^1}{W}} \leq c_\beta \, . 
\label{est.3, 3}
\end{equation}
{\cred Next, a comparison in equation \eqref{pc8} gives
\begin{equation}
\norm{\partial_t^2\ya- \fa^{(2)}}_{\vett{L^2}{H}} \leq c_\beta \, .
\label{est.3, 4}
\end{equation}
Thus, we have collected all the estimates needed to prove Theorem~\ref{th: reg1}.
In particular, note that the regularity \eqref{pc11}-\eqref{pc12} is enssured for 
the limit functions $y$ and $u$, which actually satisfy \eqref{pc13}-\eqref{pc15}.}

\bigskip
\textbf{Error equations.} We have already shown the convergence as $\alpha\searrow 
0$ for the solutions to problem $\Pab$; as we want to study the speed of
 convergence, we start writing explicitely the equations for the error. We set
\[
\yh = \ya - y \, , \qquad \uh = \ua - u \, , \qquad \xh = \xa - \xi
\]
\[
\fh = \fa - f \, , \qquad \wzh = \wza - w_0 \, , \qquad \vzh = \vza - v_0 \, , \qquad \uzh = \uza - u_0
\]
{\cred and subtract side by side the equations of the problems $\Pab$ and $\Pb$, in their strong formulations:}
\begin{equation}
\partial_t^2\yh -\beta\Delta\yh = \alpha\Delta {\cred \partial_t } \ya - \alpha\Delta\ua - \beta\Delta(1 * \uh) + \fh \qquad \textrm{a.e. in } Q
\label{eq-A_error}
\end{equation} 
\begin{equation}
\partial_t\uh - \Delta\uh + \xh + g(\ua) - g(u) = {\cred \partial_t } \yh \qquad \textrm{a.e. in } Q
\label{eq-B_error}
\end{equation}
\begin{equation}
\partial_n \yh = 0 \, , \qquad \partial_n\uh = 0 \qquad \textrm{a.e. on } \Gamma\times (0,  T)
\label{er-bd data}
\end{equation} 
\begin{equation}
\yh(0) = \wzh \, , \qquad \partial_t\yh(0) = \vzh + \uzh \, , \qquad \uh(0) = \uzh \qquad \textrm{a.e. in } \Omega \, .
\label{er-initial data}
\end{equation} 
When necessary, we also split $\fh = \fh^{(1)} + \fh^{(2)}$, where {\cred 
$\fa^{(1)}$,  $\fa^{(2)}$ (and also $f^{(1)}$,  $f^{(2)}$) fulfill \eqref{pc10},  
$\fh^{(i)} = \fa^{(i)} - f^{(i)}$ for  $i=1,2$, and 
\begin{equation}
\norm{\fh^{(1)}}_{\vett{W^{1,  1}}{V'}} +  \norm{\fh^{(2)}}_{\vett{L^1}{H}} \leq 
c\, \alpha^{1/2}. 
\end{equation}
}% fine \cred

\bigskip
\textbf{First error estimate.} Now we want to show Theorem \ref{th: error est1}. We {\cred multiply equation \eqref{eq-B_error} by $\uh$ and integrate over $Q_t$; on account of} the monotonicity of $\gamma$ and the Lipschitz-continuity of $g$, it is straightforward to get
\begin{equation}
\begin{split}
&\frac{1}{2}\nh{\uh(t)} + \int_0^t \nh{\nabla\uh(s)} ds \\
&\leq \frac{1}{2}\nh{\uzh} + c\int_0^t\nh{\uh(s)}ds + \frac{1}{2}\int_0^t\nh{\partial_t\yh(s)}ds \, .
\end{split}
\label{ts-e-1, 1}
\end{equation}
{\cred Next, we add $\beta \yh$ to both sides of \eqref{eq-A_error}, then test 
with $\partial_t\yh$ and integrate over $Q_t$. With the help of \eqref{er-initial data} and \eqref{ftc} we easily obtain}
\begin{equation}
\begin{split}
&\frac{1}{2}\nh{\partial_t\yh(t)} + \frac{\beta}{2} \norm{\yh(s)}_V^2 
\leq \frac{1}{2}\nh{\vzh + \uzh} + \frac{\beta}{2} \norm{\wzh}_V^2   \\
&\qquad {\cred {}+ \beta \norm{\wzh}_H^2 +T\beta \int_0^t \nh{\partial_t\yh(s)}ds +{}} 
\alpha\int_{Q_t}\left(\Delta{\cred \partial_t }\ya - \Delta\ua\right)
\partial_t\yh \\
&\qquad {\cred {}- \beta\int_{Q_t}\Delta (1*\uh)\, \partial_t\yh 
+ \int_0^t \dual{\fh^{(1)}(s)}{\partial_t \yh(s)} ds
+ \int_{Q_t}\fh^{(2)} \, \partial_t\yh } \, .
\end{split}
\label{ts-e-1, 2}
\end{equation}
%}% fine \cred
The term involving $\alpha$ is easily treated by the H\"older inequality and  {\cred estimates} \eqref{est.3, 2}, \eqref{est.3, 3}:
\begin{equation}
\begin{split}
\alpha\int_{Q_t}\left(\Delta{\cred \partial_t }\ya - \Delta\ua\right)\partial_t\yh \leq {\cblu \frac{\alpha}{2}} \int_0^t \alpha \nh{\Delta{\cred \partial_t }\ya(s)}ds 
+ {\cblu \frac{\alpha^2}{2}}\int_0^t \nh{\Delta\ua(s)}ds \\
{\cred {}+ \int_0^t \nh{\partial_t\yh(s)}ds 
\leq   c_\beta\, \alpha +  \int_0^t \nh{\partial_t\yh(s)}ds \, .} 
\end{split}
\label{ts-e-1, T1}
\end{equation}
{\cred Here, \eqref{pc2} has been used in the control of $\alpha^2$ by $ \alpha$. For the subsequent term, we integrate by parts in space and time and} recall \eqref{ftc}:
\begin{equation} \label{pc16}
\begin{split} 
&{\cred{}- \beta \int_{Q_t}\Delta (1*\uh)\, \partial_t\yh
\leq \beta \scal{\nabla(1*\uh)(t)}{\nabla\yh(t)} - \beta \int_{Q_t} \nabla\uh\cdot\nabla\yh } \\
&\qquad{\cred {}\leq{}} 5\beta T \int_0^t\nh{\nabla\uh(s)} ds + \frac{\beta}{8}\nh{\nabla\yh(t)} + \frac{\beta}{4T}\int_0^t \nh{\nabla\yh(s)} ds \, .
\end{split}
\end{equation}
{\cblu The next term in \eqref{ts-e-1, 2} is treated similarly as in \eqref{term-f1}. We infer that 
\begin{equation}\label{pc17}
\begin{split}
\int_0^t \dual{\fh^{(1)}(s)}{\partial_t \yh(s)} ds
\leq c_\beta \norm{\fh^{(1)}}^2_{\vett{L^\infty}{V'}} 
+ \frac{\beta}{8}\nh{\nabla\yh(t)} +\frac{\beta}{4}\nh{\wzh}\\
{}+ \frac{\beta T}{4}\int_0^t \nh{\partial_t\yh(s)} ds
+ \frac{1}{2}\norm{\wzh}^2_V 
+ \int_0^t\norm{\partial_t \fh^{(1)}(s)}_{V'}\norm{\yh(s)}_V ds 
 \, ,
\end{split}
\end{equation}
}% fine \cblu
while the other contribution of the source term is easily estimated by the H\"older inequality:
\begin{equation} \label{pc18}
\int_{Q_t} \fh^{(2)}\partial_t\yh \leq \int_0^t \norm{\fh^{(2)}(s)}_H\norm{\partial_t\yh(s)}_H ds \, .
\end{equation}

Now we multiply \eqref{ts-e-1, 1} {\cred by $(1 + 5\beta T) $ and add the resulting inequality to \eqref{ts-e-1, 2}, side by side. In view of \eqref{ts-e-1, T1}--\eqref{pc18} and thanks to \eqref{f rate est1}--\eqref{d rate est1}, 
we deduce that
\begin{equation}
\begin{split}
&\frac12 \nh{\uh(t)} + \int_0^t \nh{\nabla\uh(s)} ds + \frac{1}{2}\nh{\partial_t\yh(t)} + \frac{\beta}{4} \norm{\yh(t)}_V^2 \\
&\leq c_\beta\, \alpha +  c\int_0^t\nh{\uh(s)}ds + c\int_0^t\nh{\partial_t\yh(s)}ds + c\int_0^t \nh{\nabla\yh(s)}ds \\
&\qquad + \int_0^t\norm{\partial_t \fh^{(1)}(s)}_{V'}\norm{\yh(s)}_V ds + \int_0^t \norm{\fh^{(2)}(s)}_H\norm{\partial_t\yh(s)}_H ds  \, ;
\end{split}
\label{ts-e-1, 3}
\end{equation}
}% fine \cred 
hence, the generalised Gronwall lemma (see{\cred , e.g.,}  \cite{baiocchi}) and \eqref{f rate est1} again entail the thesis of Theorem \ref{th: error est1}.

\bigskip
\textbf{A regularity result on $u$.} {\cred Let us formally describe how  Theorem~\ref{th: reg2} can be proved: you differentiate the equation in \eqref{pc8} with respect to time and then test by $\partial_t \ya$, integrating over $Q_t$ with the help of the further initial condition for $\partial_t \ya (0)$. Such a condition can be still inferred from \eqref{pc8} and, thanks to \eqref{dn reg2}, the initial values $\partial_t \ya (0)$ turn out to be bounded in $H$. Moreover,
you exploit the facts that 
$$ \int_{Q_t} \partial_t \xa \, \partial_t \ya \geq 0$$
and that $ \partial_t^2 \ya $ stays uniformly bounded in $\vett{L^1}{H}$ 
(cf.~\eqref{est.3, 4} and \eqref{f reg1}).

However, for a rigorous proof we refer to \cite[Section~7]{cc1}, where the sketched estimate is carried out properly, in the setting provided by the Faedo-Galerkin approximation and the Yosida regularisation of the graph $\gamma$. Since the estimate is independent of $\alpha$, for the proof, as a by-product we have 
the desired inequality
\begin{equation}
\norm{\ua}_{\vett{W^{1,  \infty}}{H} \cap \vett{H^1}{V} \cap \vett{L^\infty}{W}} \leq c_\beta \, ,
\label{reg2}
\end{equation}
which pass to the limit as $\alpha\searrow 0$ via compactness argument.}

\section{The case of a smooth potential}
\label{sm-po}
In what follows, we assume that $\gamma$ is a single-valued, non decreasing, 
locally Lipschitz-continuous function, defined on the whole {\cred of $\R$ (see 
\eqref{gamma} and \eqref{gamma smooth}). {\cblu We have already pointed out that this assumption is restrictive but not unrealistic, as the example of the `double-well' potential \eqref{double well} shows.}

Since the hypotheses we consider on the data are strong enough to guarantee \eqref{reg2}, and consequently the estimate 
\begin{equation}
\norm{\ua}_{L^{\infty} (Q) } \leq c_\beta \, ,
\nonumber
\end{equation}
we have that $\xa = \gamma(\ua) $ is bounded independently of $\alpha$
in $ L^\infty(Q)$. 
Furthermore, by the local Lipschitz-continuity of $\gamma$, without loss of 
generality we can set  $\xa = 0$ in \eqref{eq-B_alpha} (and \eqref{pc8}), 
$\xi = 0$  in \eqref{eq-B} (and \eqref{pc14}), $\xh = 0$
in \eqref{eq-B_error} by modifying $g$ in order to take account of $\gamma$ too.

\bigskip
\textbf{Further regularity on $u$.} As we are interested in showing Theorem \ref{th: reg2.5}, we consider equation \eqref{eq-B_alpha} and derive it, with respect to time. We get a first-order equation for a new unknown $z_\alpha = \partial_t\ua$, which is complemented by Neumann homogeneus boundary conditions 
and the suitable initial condition, obtained by setting $t = 0$ in \eqref{pc8}:
\[
\partial_t z_\alpha - \Delta z_\alpha = -\partial_t (g(\ua)) + \partial_t^2 \ya \quad \textrm{a.e. in } Q \, ,  
\]
\[
\partial_n z_\alpha = 0 \quad \textrm{a.e. in } \Gamma\times (0,  T)
\]
\[
z_\alpha(0) = \Delta\uza -  g(\uza) + \vza + \uza \quad \textrm{a.e. in } \Omega \, .
\]
We notice that the terms $\partial_t (g(\ua))$ and $ \partial_t^2 \ya$ are bounded 
in $ L^2(Q)$ due to \eqref{g}, \eqref{est.3, 2} and 
\eqref{est.3, 4}, \eqref{d reg2.5}. Moreover, \eqref{dn reg1} and \eqref{d reg2.5}
imply that
\[
\Delta\uza - g(\uza) + \vza + \uza \    \hbox{ is bounded in } \,  V 
\]
independently of $\alpha$. Thus, we can apply Lemma~\ref{lemma: parabolic} 
and, on account of \eqref{reg2}, the inequality \eqref{parabolic} leads the estimate 
\begin{equation}
\norm{\ua}_{\vett{ H^2}{H} \cap \vett{ W^{1,  \infty}  }{V} \cap \vett{H^1}{W}} \leq c_\beta \, ,
\label{pc19}
\end{equation}
whence the thesis of Theorem \ref{th: reg2.5} follows.

\bigskip
\textbf{Further regularity on $\partial_t y$.} Here, we want to prove Theorem~\ref{th: reg3}, whose assumptions are now supposed to be in force. 
Recalling \eqref{pc7}, \eqref{pc9} and \eqref{i-c-alpha}, it turns out that $v_\alpha = \partial_t \ya $ formally satisfies 
\begin{equation}\label{pc20}
\partial_t^2 v_\alpha - \alpha\Delta \partial_t v_\alpha - \beta\Delta v_\alpha = -\alpha\Delta \partial_t \ua -\beta\Delta \ua + \partial_t \fa  \quad \textrm{ a.e. in } Q
\end{equation}
\begin{equation}\label{pc21}
\partial_n v_\alpha = 0 \quad \textrm{ a.e. on } \Gamma\times (0, T) 
\end{equation}
\begin{equation}\label{pc22}
v_\alpha ( 0)= \vza + \uza, \quad  \partial_t  v_\alpha (0 )= 
\alpha \Delta \vza +\beta \Delta \wza + \fa(0)  \quad \textrm{ a.e. in } \Omega .
\end{equation}
Due to \eqref{pc19} and \eqref{f reg3}, the right-hand side in \eqref{pc20}
has exactly the same regularity \emph{here} as its counterpart in \eqref{pc7} \emph{within the framework of Theorem~\ref{th: reg1}}. The same correspondence holds for the initial data in \eqref{pc22} and the related ones in \eqref{i-c-alpha}: indeed, 
$\vza + \uza$ is in $W$ and $ \alpha \Delta \vza +\beta \Delta \wza + \fa(0)$ lies in $V$ due to \eqref{d reg2.5} and \eqref{d reg3}. Then, we have that (cf.  \eqref{pc5})  
\begin{equation}
v_\alpha \in \vett{W^{2,1}}{H} \cap \vett{W^{1, \infty}}{V} \cap \vett{H^1}{W}
%\nonumber
\nonumber
\end{equation}
and $v_\alpha$ actually satisfies \eqref{pc20}--\eqref{pc22}. Then, we are allowed
to test \eqref{pc20} by $-\Delta\partial_t v_\alpha $ and repeat the computations 
developed in the Fourth a priori estimate. Thanks to  \eqref{d reg2.5} and 
\eqref{dn reg3}, proceeding in this way we obtain (cf.~\eqref{pc23})
\[
\begin{split}
\frac{1}{2} \nh{ \nabla \partial_t v_\alpha(t)} 
+ \frac{\alpha}{2} \int_0^t \nh{\Delta\partial_t v_\alpha(s)} ds 
+ \frac{\beta}{4} \nh{\Delta v_\alpha(t)} 
\leq  c_\beta + c \int_0^t \nh{\Delta  v_\alpha(s)} ds \\
+ \int_0^t\norm{\partial_t^2 \fa^{(1)}(s) }_H \norm{\Delta v_\alpha(s)}_H ds 
+ \int_0^t \norm{\nabla \partial_t  \fa^{(2)}(s)}_H 
\norm{\nabla\partial_t v_\alpha(s)}_H ds \, .
\end{split}
\]
By means of the generalised Gronwall lemma, stated in {\cred \cite[Teorema~2.1, p.~245]{baiocchi}}, recalling \eqref{f reg3}, \eqref{dn reg3} and $v_\alpha = \partial_t \ya$, we can conclude
\begin{equation}
\norm{\partial_t \ya}_{\vett{W^{1,  \infty}}{V} \cap \vett{L^\infty}{W}} + \sqrt{\alpha}\norm{\partial_t\ya}_{\vett{H^1}{W}} \leq c_\beta \, ,
\label{reg3}
\end{equation}
which entails \eqref{reg 3}.

\bigskip
\textbf{Second error estimate.} We begin with revisiting the first estimate 
on the convergence error. Assume all the hypotheses of Theorem~\ref{th: error est2}, and reconsider the proof of Theorem~\ref{th: error est1} we have given in 
Section~\ref{sec: strong}. We leave unchanged all the estimates we have already established, except the inequalities~\eqref{ts-e-1, T1}: indeed, from \eqref{est.3, 2} and \eqref{reg3} it follows that
\[
\begin{split}
\alpha\int_{Q_t}\left(\Delta \partial_t \ya - \Delta\ua\right)\partial_t\yh \leq {\cblu\frac{\alpha^2}{2}}\int_0^t \nh{\Delta\partial_t\ya(s)}ds  + {\cblu\frac{\alpha^2}{2}} \int_0^t \nh{\Delta\ua(s)}ds \\
{}+ \int_0^t \nh{\partial_t\yh(s)}ds 
\leq  c_\beta\, \alpha^2 + \int_0^t \nh{\partial_t\yh(s)}ds  \, .
\end{split}
\]
By this new estimate and the assumptions~\eqref{f rate est1-a}--\eqref{d rate est1-a}, inequality~\eqref{ts-e-1, 3} now becomes
\[
\begin{split}
&\frac12 \nh{\uh(t)} + \int_0^t \nh{\nabla\uh(s)} ds + \frac{1}{2}\nh{\partial_t\yh(t)} + \frac{\beta}{4} \norm{\yh(t)}_V^2 \\
&\leq c_\beta\, \alpha^2 +  c\int_0^t\nh{\uh(s)}ds + c\int_0^t\nh{\partial_t\yh(s)}ds + c\int_0^t \nh{\nabla\yh(s)}ds \\
&\qquad + \int_0^t\norm{\partial_t \fh^{(1)}(s)}_{V'}\norm{\yh(s)}_V ds + \int_0^t \norm{\fh^{(2)}(s)}_H\norm{\partial_t\yh(s)}_H ds  \, ;
\end{split}
%\nonumber
\]
so that the generalised Gronwall lemma {\cred (cf., e.g., \cite{baiocchi})} and \eqref{f rate est1-a} entail
\begin{equation}
\norm{\yh}_{\vett{W^{1,  \infty}}{H} \cap \vett{L^\infty}{V}} + \norm{\uh}_{\vett{L^\infty}{H} \cap \vett{L^2}{V}} \leq c_\beta\, \alpha \, .
\label{er-est1, alpha}
\end{equation}
Next, we consider the equation \eqref{eq-B_error}, where we have set $\xh = 0$, and apply Lemma \ref{lemma: parabolic}; by \eqref{d rate est1-a}, the Lipschitz-continuity of $g$ and \eqref{er-est1, alpha} we obtain
\begin{equation}
\norm{\uh}_{\vett{H^1}{H} \cap \vett{L^\infty}{V} \cap \vett{L^2}{W}} \leq c\left\{ \norm{\uzh}_V + \norm{\uh}_{L^2(Q)} + \norm{\partial_t\yh}_{L^2(Q)} \right\} \leq c_\beta\, \alpha \, .
\label{er-est 2, 1}
\end{equation}
Then \eqref{pc31} is entirely proved.  In order to complete 
the proof of Theorem \ref{th: error est2}, we test equation 
\eqref{eq-A_error} with $-\Delta\partial_t\ya$ (which belongs to $L^2(Q)$ because of 
the estimate \eqref{reg3}) and integrate over $Q_t$:
\begin{equation}
\begin{split}
\frac{1}{2}\nh{\nabla\partial_t\yh(t)} + \frac{\beta}{2}\nh{\Delta\yh(t)} \leq \frac{1}{2}\nh{\nabla \left( \vzh + \uzh \right) } + \frac{\beta}{2}\nh{\Delta\wzh} \\
+ \alpha\int_{Q_t} \left(\Delta \ua - \Delta\partial_t\ya\right)\, \Delta\partial_t\yh + \beta\int_{Q_t}\Delta(1 * \uh)\, 
\Delta\partial_t\yh - \int_{Q_t}\fh \, \Delta\partial_t\yh \, .
\end{split}
\label{ts-error 2, 1}
\end{equation}
We deal with the integral involving $\ua$ and $\ya$, which is delicate to treat. By an application of the H\"older inequality, \eqref{reg2}, and  \eqref{reg3} -- which holds for $y$ as well -- we have 
\[
\begin{split}
&\alpha\int_{Q_t} \left(\Delta\ua - \Delta\partial_t\ya\right)\, \Delta\partial_t\yh \\ &{}\leq {\cblu \alpha\left( \norm{\Delta\ua}_{L^2(Q)} + \norm{\Delta\partial_t\ya}_{L^2(Q)} \right) \left(\norm{\Delta\partial_t\ya}_{L^2(Q)} + \norm{\Delta\partial_t y}_{L^2(Q)} \right)} \leq c_\beta\, \alpha \, .
\end{split}
\]
The subsequent term is estimated by integrating by parts, with respect to time, and taking advantage of \eqref{ftc}:
\[
\begin{split}
\beta\int_{Q_t}\Delta(1 * \uh)\, \Delta\partial_t\yh \leq 2\beta \nh{\Delta(1 * \uh)(t)} + \frac{\beta}{8}\nh{\Delta\yh(t)} - \beta\int_{Q_t}\Delta\uh \, \Delta\yh \\
\leq \beta \left(4 T + 1\right) \int_0^t \nh{\Delta\uh(s)} ds + \frac{\beta}{8}\nh{\Delta\yh(t)} + \frac{\beta}{4} \int_0^t \nh{\Delta\yh(s)} ds \, .
\end{split}
\]
Now we split $\fh = \fh^{(1)} + \fh^{(2)}$, as usual, and deal on the respective terms separately, integrating by parts the former in time and the latter in space. We obtain
\[
\begin{split}
\int_{Q_t}\fh^{(1)}\, \Delta\partial_t\yh \leq \frac{2}{\beta}\nh{\fh^{(1)}(t)} + \frac{\beta}{8}\nh{\Delta\yh(t)} + \frac{1}{2}\nh{\fh^{(1)}(0)} \\
{} + \frac{1}{2}\nh{\Delta\wza} + \int_0^t\norm{\fh^{(1)}(s)}_H\norm{\Delta\yh(s)}_H ds
\end{split}
\]
and
\[
\begin{split}
\int_{Q_t}\fh^{(2)}\, \Delta\partial_t\yh \leq \int_0^t\norm{\nabla\fh^{(2)}(s)}_H\norm{\nabla\partial_t \yh(s)}_H ds \, .
\end{split}
\]
According to all these estimates, and taking \eqref{d rate est2}, \eqref{er-est 2, 1} and \eqref{pc2} into account, the inequality \eqref{ts-error 2, 1}
transforms into 
\[
\begin{split}
&\frac{1}{2}\nh{\nabla\partial_t\yh(t)} + \frac{\beta}{4}\nh{\Delta\yh(t)} \leq  c_\beta\, \alpha + \frac{\beta}{2} \int_0^t \nh{\Delta\yh(s)} ds \\
&\qquad{}+ \int_0^t\norm{\fh^{(1)}(s)}_H\norm{\Delta\yh(s)}_H ds + \int_0^t\norm{\nabla\fh^{(2)}(s)}_H\norm{\nabla\partial_t\yh(s)}_H ds  \, .
\end{split}
\]
In view of this inequality and \eqref{d rate est2}, by 
the generalised Gronwall lemma {\cred (cf., e.g., \cite{baiocchi})} we
deduce \eqref{pc32} and therefore conclude the proof of Theorem~\ref{th: error est2}.
}% fine \cred

%%%%%%%%%%%%%%%%%%%%%%%%%%%%%%%%%
%% bibliography
%%%%%%%%%%%%%%%%%%%%%%%%%%%%%%%%%

%\vspace{3truemm}

\end{document}